Bayes´ Theorem and the cMPE-Method: Differences, Complementarities
Inconsistencies and the Importance of Discerning between Weight-Bearing and Extensional Evidence.
By R. Gottlob

I. INTRODUCTION

Bayes´ Theorem (BayT) is regarded as the only possibility of quantitatively increasing probabilities because of new relevant experience. Recently, a second method was published, the cMPE-Method, permitting the adding of probabilities, however, taking into account their non-linearity(R. Gottlob, 2000, 2003 2003). The method is based on the **m**ultiplication of the **p**robabilities of **e**rror [$P(\neg X) = 1 - P(X)$]. The product of the probabilities of error is subtracted from 1. A precondition is `*semantic*´ independence. Analogously, the DPE-Method permits subtractions of probabilities, in spite of their non-linearity. These methods do not compete with BayT, but rather, they are complementary.

In Section II BayT is discussed, some advantages are outlined and some shortcomings are mentioned. In Section III, the MPE-Method, the cMPE-Method and the DPE-Method are described and in Section IV the differences between the BayT and the MPE-Methods are outlined. Section V is devoted to the fundamental differences between extensional and weight-bearing[1] evidence. Some unsolved problems with BayT will be described.

Some notations and formulas taken for granted:
The probability of the coincidence of two stochastically independent data A and B is

$$P(A \& B) = P(A) \times P(B) \qquad (1a)$$

In case of stochastic dependence:

$$P(A \& B) = P(A|B) \times P(B) \text{ or } P(A) \times P(B|A) \qquad (1b)$$

Instead of A and B we may use H and E for hypothesis and evidence.

The complement of $P(X)$, $1 - P(X)$ may be called also as the probability of error of X and the notation is $P(X^c)$. Additions or subtractions of probabilities, however, taking into account their non-linearity, are notated by the superscripts $^+$ and $^-$. In all probabilistic reasoning our background knowledge K must be taken into account. For the sake of clarity this is done, however tacitly, in the following considerations.

II. BAYES´ THEOREM

is derived from (1b).

$$P(H|E) = P(H) \times P(E|H) / P(H) \times P(E|H) + P(\neg H) \times P(E|\neg H) \qquad (2a)$$

The denominator of (2a) contains the total probability of E . More generally, BayT with the total probability reads:

---

[1] J.M. Keynes(1957, ch. VI) discusses the weight of arguments . Hereby, Keynes comprises the absolute amount of relevant knowledge and of relevant ignorance, respectively. According to Keynes, New evidence will sometimes decrease the probability of an argument, but it will always increase the `weight´ . Here, I shall speak only of the ´weight-bearing´ of evidence in order to indicate evidence that only *increases* probabilities.

$$P(A_k | B) = \frac{P(B | A_k) \times P(A_k)}{\sum_{i=1}^{n} P(B | A_i) \times P(A_i)} \quad (k = 1, \ldots, n) \quad (2b)$$

Fig: 1: The BayT formula with total probability.

According to R. v. Mises, all probabilities may be derived from relative frequencies, represented by the quotient of the number of favourable or relevant cases, divided by the number of possible cases. Also BayT may be derived from relative frequency: Both, the numerator and the denominator are formed by products, therefore, their absolute values are smaller than P(A) or P(B). However, the quotient is assumed to equal P(A|B) or, for hypothesis H and evidence E, P(H|E). The more extended the probability in the denominator, the lower P(H|E). It is an advantage of BayT that the important cases of E| ¬H are also allowed for.
Fig. 2 will schematically show some possible applications of BayT:

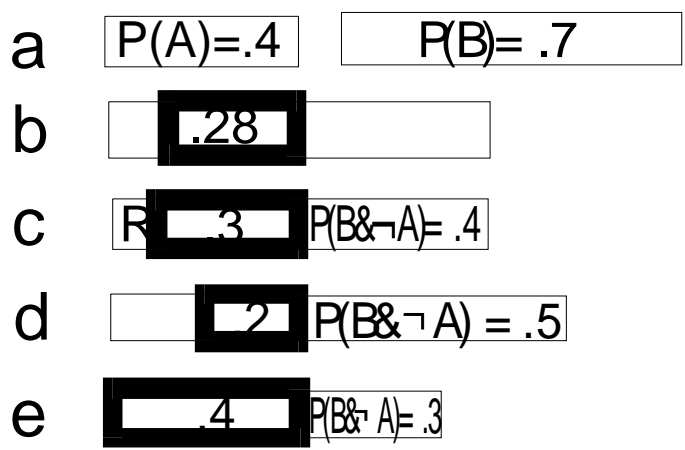

Fig. 2: Some applications of BayT. P(A & B) or P(A) x P(B|A) are indicated by boldfaced outlines. **a:** Premises (.4, .7) **b:** Overlapping by sheer chance (stochastic independence). **c:** The overlapping domain is .3, P(A|B) = .43, (as the quotient .3 : .7). However note that the domain R has no share in this quotient, because P(A)x P(B|A) < P(A). **d**: P(A) x P(B|A) is small, P(B & ¬A)is large, therefore P(A|B) = .2 / .7 = .29 and <P(A), <P(B). **e:** Hypothesis A implies evidence B.

An ambiguity of BayT is shown in Fig. 2**c**. What is correctly computed is a division of P(A) x P(B|A) by the total P(B) = .7. However, P(A) x P(B|A) < P(A) and there remains a rest `R´ unconsidered. This calculation is correct if we wish to know the actual chance of coincidence of the total P(B) with P(A) x P(B|A). However, in most instances, what we wish to know is the increase of P(A) under the influence of P(B) and this can only be .43 $^+$ R . A simple addition, P( R) + .43 is not permitted, because sometimes, P(B & ¬A) may be but slightly larger than P(A) x P(B|A) and then

the quotient may be almost 1 or close to 1 and by simple addition of R, probability 1 may be exceeded. The correct addition of such nonlinear probabilities will be described below.
As will be detailed in section V, severe flaws may arise, if for implication of B by A (Fig. 2**e**) the differences between extensional Es and weight-bearing Es are disregarded.

III: THE MPE-; THE cMPE- AND THE DPE-METHOD.
C. Glymour (1980, ch. III), under the heading `Why I am not a Bayesian´ pointed out that BayT was deemed the only method available to increase the probability of a hypothesis due to some new evidence. There might be other procedures, but "*none of these arguments is decisive against the Bayesian scheme of things, nor should they be, for in important respects that scheme is undoubtedly correct. But taken together, I think they do at least strongly suggest that there must be relations between evidence and hypotheses that are important to scientific argument and to confirmation but to which the Bayesian scheme has not yet penetrated*". (p. 93).
Such an argument may be given by the MPE- and by thc cMPE-Method.
A: *The MPE-Method (**M**ultiplication of the **P**robabilities of **E**rror)*
consists in the multiplication of the probabilities of error, whereby the product will become smaller than the single factors.

$$P(\neg A) \,\&\, P(\neg B) \,\&\, \ldots \,\&\, P(\neg N) = P(\neg A) \times P(\neg B) \times \ldots \times P(\neg N) \qquad (3)$$

This operation may be used consciously or unconsciously for many types of cognition, e.g. for singular facts. (4) is an analogue to 1, however, it cannot be based on stochastic, but only on *semantic* independence (see below).
*Example:* `This object in my hand is an apple´ (i) I recognize it by its visual appearance. Until now, I certainly saw more than 2000 apples; I remember only two times to have seen a dummy, so I guess that the probability of error will be at most 1/1000. (ii) The purely tactile impression (weight, surface, humidity) could lead us astray at most in one case of 100 and (iii) the very specific smell at most in one case of 1000. By multiplication $10^{-3} \times 10^{-2} \times 10^{-3}$ we arrive at an at-most value of $10^{-7}$ as probability of error and this low value prevented us to bite into a pebble, if we wanted to eat an apple.
The role of the MPE-Method for cognition of simple facts and especially of Gestalt cognition will be extensively discussed elsewhere. In this paper we shall be mainly concerned with
B: *The cMPE-Method*
As all probabilities are fractions <1, a multiplication like in (1a) or in (1b) will lead to decreased probabilities. However, if two or more probabilities support one another, an increased probability must result. This may be accomplished by multiplying the complements of the probabilities and by subtracting the product from 1.

$$P(A) \,^{+}\, P(B) = 1 \,[1-P(A)] \times [1 \, P(B)] \qquad (4a)$$

Or

$$P(A) \,^{+}\, P(B) = 1 \, P(A^c) \times P(B^c) \qquad (4b)$$

Or

$$P(A) \,^{+}\, P(B) \,^{+}\, \ldots \,^{+}\, P(N) = 1 \, P(A^c) \times P(B^c) \times \ldots \times P(N^c) \qquad (4c)$$

Or, after carrying out the multiplication in (5a)

$$P(A) \,^{+}\, P(B) = P(A) + P(B) \, P(A) \times P(B) \qquad (4d)$$

The right-hand side of (4d) resembles formulas for the union or for at least one event being true in a series of stochastically independent events. However, in (4a - d) the events are stochastically dependent, and on the left-hand side an addition, however, taking into account the non-linearity of probabilities is indicated (see also Fig. 4). No probability >1 may arise. (4a - d) are formulas, expressing what W. Whewell (1840, pp. 230, 232) called

Consilience of Induction : *"According to cases in which induction from classes of facts altogether different have jumped together, belong to the best established theories which the history of science contains. No example can be pointed out, in the whole history of science, so far as I am aware, in which this Consilience of induction has given testimony in favour of an hypothesis afterwards discovered to be false"* . There, classes of facts, altogether different cannot mean stochastic independence but only,

*C: Semantic variety* ( Greek: σημάιναι indicate) As the apple-example demonstrated, our cognition is based on different ways of perception or on different properties of the objects. Semantic variety permits multiplications of the probabilities of error. Besides different sense organs, different scientific methods may be of importance, and in law courts, testimonies of independent witnesses supporting one another or being supported by circumstantial evidence are deemed highly reliable. However, stochastic independence as prerequisite for a descending operation, such as (1a) cannot be assumed for ascending operations (such as for BayT and for the cMPE-Method), where A supports B etc. However, there must be no overlapping of the pieces of evidence besides their jointly cionfirming a hypothesis.

Variety of evidence plays an important role, or even more important than sheer amount of evidence (J. Earman, 1992, p. 77 and passim, C. Glymour, 1980, p.139,ff.). According to G. Boole (1953, p. 155):

*"It is obviously supposed ( .) that the probabilities of the simple events x, y, z, are given by independent observations. This is, I apprehend, what is really meant by events being spoken as independent (according to the ordinary acception of that term), the knowledge that they are not so can only be derived from experience in which they are mutually involved, not from observations upon them as simple and unconnected, and hence, if our knowledge is derived from experience, the independence of events is only another name for the independence of our observation of them".* And on page 157 we read:

*"The events, whether simple or compound, whose probabilities are given by observation, are to be regarded as independent of any but a logical connection".*

C. Glymour (1980,p. 140) notes, that different hypotheses may make one way of testing relevantly different from others. *"Part of what makes one piece of evidence relevantly different from another piece of evidence is that some test is possible from the first that is not possible from the second, or that in the two cases there is some difference in the precision of computed values of theoretical quantities".*

With these statements, Glymour comes rather close to what here is called semantic independence.

D: *Derivation*

It is generally agreed that $P(A) \times P(B) + P(\neg A) \times P(B) + P(A) \times P(\neg B) + P(\neg A) \times P(\neg B) = 1$.
Therefore:
$$P(A) \times P(B) + P(\neg A) \times P(B) + P(A) \times P(\neg B) = 1 - P(\neg A) \times P(\neg B) \qquad (5)$$
The left-hand side of (5) contains all positive values of P(A) or P(B), multiplied by a factor that prevents arising probabilities >1. For checking we assume that P(A) = .3 and P(B) = .8.
Then we arrive for the left-hand side of (5) at .24 + .56 + .06 = .86 and for the right-hand side we arrive at 1 - .14 = .86. For more about derivation of (5a d) and for relations to percentage operations see Gottlob (2003).

E: *Complements of the results:*
The complements of the results of 5(a d) may be calculated e.g., by
$1-[1-P(A^c) \times P(B^c)] = P(A^c) \times P(B^c)$ . It must be decided by our background knowledge, whether we apply in our calculations the methods 4 (a - d) or whether we just multiply the complements $(P(A^c) \times P(B^c))$. L.J. Cohen (1977, p. 100 ff.), in order to falsify the percentual calculations of P.O. Ekelöf´s (1964), assumes the testimonies of two witnesses, both reliable only to 25 %. Both had claimed that a criminal was a male (Witness A: Because he had long hair; witness B was known to be an ardent feminist). These testimonies result with 1 - .75 x .75 = .44 % in favour of maleness. Because the probability of maleness must be the complement of femaleness, we may calculate also: 1 - .25 x .25 = 94 for femaleness. This leads to the relation 44% maleness and 94 % femaleness, whereas the numbers should be complementary. This absurd relation arouse because right at the beginning 25 % (P = .25) for maleness had to be regarded as the complement of the more believable 75 % for femaleness. By simple multiplication of the complements (.25 x .25) we arrive at .06, and this is the real complement of .94 for femaleness. We may conclude that in binary cases, the larger probabilities must be calculated according to (4a - d), whereas the lower probabilities must be regarded as complements and must be calculated by simple multiplication. (See Gottlob 3003).

F: *Some ways of application*

Fig. 3 shows some ways of application of the cMPE-Method schematically.

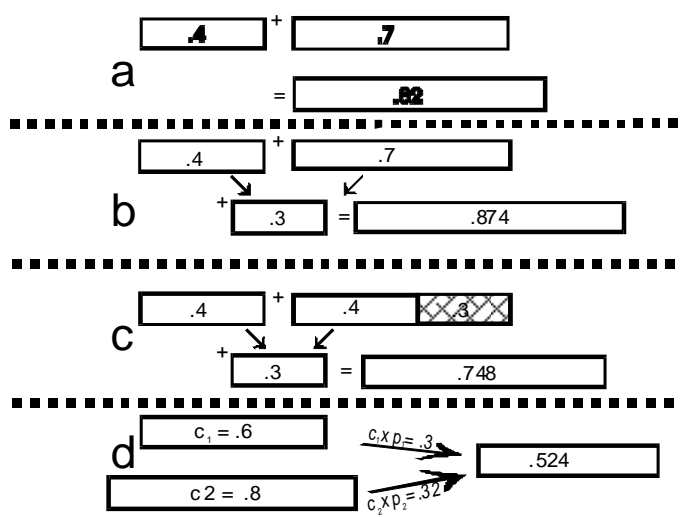

Fig. 3 **a**: Event A [P(A) = .4] supports event B [P(B)=.7] totally, .4 $^+$ .7 = .82 result.
**b:** Event A (.4) and event B (.7) jointly support event C (.3) .4 $^+$ .7 $^+$ .3 = .874 result.
**c:** A similar situation as in **b**, however, event B only partially (with probability .4) supports A and C . **d**: The a priori probability of A, $c_1$ = .6. The support of C by A, $p_1$ = .5. Thus A, according to (1 b), confers to C (that has no a priori probability) probability .6 x .5 =.3 The a priori support of B, c2 = .8. The support of C by B, p2 = .4 . Thus B confers to C the probability .8 x .4 = .32 and the resulting probability of C will be .3 $^+$ .32 = 1 - .7 x .68 = .524.

In Fig. 3 the examples **a - c** are self-evident. Example **d**, however, deserves some comments. J.M. Keynes (1952, ch. XVII) criticises the rather involved formula of `G. Boole´s. Challenge Problem´. Keynes´ own formula appears to be simpler, but at the expense of introducing new symbols like u, y or z. Our formula for P( C), an event that has no probability of its own, but obtains some probability $p_1$ and $p_2$ from the events A or B, respectively, that have the a priori probabilities $c_1$ and $c_2$ reads:

$$P(C) = P(c_1) \times P(p_1|c_1) \; ^+ \; P(c_2) \times P(p_2|c_2) \tag{6}$$

What follows, is a simple example of an application of the cMPE-Method:

*Example:* In a certain region the occurrence of thunderstorms within one weak in July is .6 If the average day temperature is exceeded by $2^oC$ the probability of thunderstorms increases by .4 and if the average air humidity is increased by 20 % the probability of thunderstorms increases by 50 %. What is the probability of thunderstorms if both, the day temperature and the air humidity increase as mentioned? A simple addition cannot be performed, because adding .4 and .5 to .6 would result in 1.5. By the cMPE-Method we arrive at 1 - .4 x .6 x .5 = 88.

G: *The DPE-Method*

Non-linear subtractions of probabilities may be performed by the related DPE Method:

$$P(A) \; ^- \; P(B) = 1 - P(A^c) / P(B^c) \tag{7a}$$

For a multitude of subtrahends we obtain:

$$P(A) \; ^- \; P(B) \; ^- \; .... \; ^- \; P(N) = 1 \; P(A^c) / P(B^c) \times .... \times P(N^c) \tag{7b}$$

*Example:* A certain disease may be cured by a new drug in 60 % (P = .6), however in 20 % (P = .2), some untoward side effects are observed. A conventional drug permitted cures in 50 %, but in 10 % side effects occurred. The difference between the cures and the side effects may be called the `therapeutic window´. Which of the treatments has a wider therapeutic window?

We calculate: $.5 \; ^- \; .1 = 1 - .5 / .9 = .444$ for the conventional drug,
but $.6 \; ^- \; .2 = 1 - .4 / .8 = .5$ for the new treatment.

H: *The role of non-linearity.*

The above example shows that due to the non-linearity, higher probabilities by subtracting certain values are changed to a lesser degree than low probabilities. This is demonstrated for additions in Fig. 4:

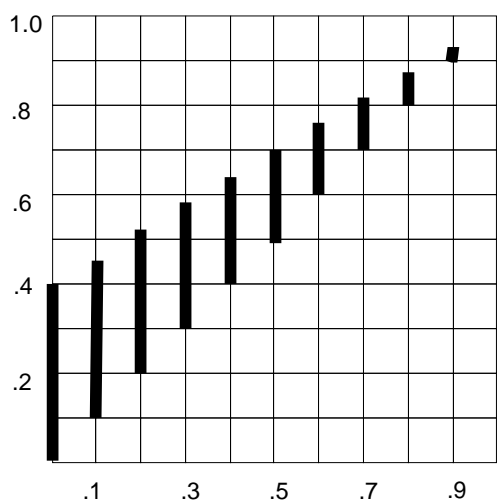

Fig. 4: To the values of the abscissa, .4 is added by the cMPE-Method. Note that due to the non-linearity high probabilities are changed least. The total added value of .4 is preserved only for the initial value 0.

Due to the non-linearity no probabilities >1 may result and by adding 1.0 to any probability, the sum will be 1.0. Rather high probabilities may be added to high probabilities with a seemingly insignificant effect: If we add .999 to a probability of .99, we arrive at .999 99 .However, subtracting .999 from .999 99 results in 1 - .00001 / .001 = .99 according to the DPE-Method, whereby the necessary great input for causing minor changes in large probabilities is proved..
By taking into account the non-linearity, the Puzzle of old evidence (e.g., C. Glymour, 1980,) may be solved.
The puzzle consists in the fact that highly proved new evidence (say, the General Theory of Relativity) cannot be further proved significantly by some important old evidence, that was not considered while the theory was elaborated (e.g. the advance of the perihelion of Mercury). The explanation is given above: Adding non-linearly a large probability to another large probability causes but insignificant linear increases. Subtracting the increment again results in high probabilities. In other words, very small linear changes in high probabilities may reflect considerable changes if calculated by the appropriate non-linear method.
I: *A comparison with Laplace´s Method and combined applications.*
In this subsection we shall consider the use of Laplace´s Method for preventing overly arbitrary estimations of prior probabilities (in BayT) or of any probability that serves as premise for the cMPE-Method. True, several authors, such as J.M. Keynes (1957, p. 377 ff and passim) put forward severe objections against this method and indeed, for very low numbers of observation, the method is unreliable. Thus it makes little sense to assume after one observation (say, of a white ball drawn from an urn) that the probability for the next ball to be white will be 2/3 or after observation of two white balls, to assume that the probability for the next ball to be white is ¾ . However, the reliability of the rule : P = (n +1) / (n + 2) increases with increasing numbers of observation. (For a derivation, see B. Russell, 1948, p. 367). Now even authors that have reservation against the method of relative frequency will admit that this method is still preferable to unfounded hunches. As our calculations are mostly devoted to processing probabilities obtained from large numbers of observation, frequently without any unexplained counterexample observed, I propose to apply wherever possible a combined method for the cMPE-Method as well as for assessing the prior of the BayT.

*Example:* Assume that at some time it was found that among birds in Europe the first 50 observed animals were lung breathers. The probability that the next bird encountered would be a lung breather was 51 / 52 = .98 .After 100 such conform observation this probability could be increased to .99. However, if the first 50 animals were European birds and the second 50 birds observed came from tropic regions and were fundamentally different in several respects, but all were lung breathers, one could have calculated according to the cMPE-Method $P(A)^+ P(B) = 1 - .02 \times .02 = .9996$. At first blush this difference: 99 for homogeneous observations and .9996 for observations of two different groups appears as not too significant. However, by the DPE-Method we learn, that $.9996^- .99 = 1 - .0004 / .01 = .96$.

For a comparison of the results of Laplace´s Method and the cMPE-Method see Table 1:

Table 1

| 1    | 2       | 3       | 4       | 5                  |
|------|---------|---------|---------|--------------------|
|      | Laplace | Laplace | cMPE    | cMPE 2A$^-$        |
| A    | A       | 2A      | 2A      | Laplace 2A         |
| 5    | .857..  | .917    | .979    | .253               |
| 10   | .917    | .955    | .993    | .286               |
| 50   | .98     | .990    | .996    | .96                |
| 100  | .99     | .995    | .999 9  | .98                |
| 1000 | .999    | .999 5  | .999 999| .998               |

Table 1: A comparison of the Laplacean Method for homogeneous observations with the cMPE-Method for two semantically diverse collectives. A: Numbers of observation. The increasing numbers in Column 5 indicate a marked superiority of the cMPE-Method in detecting probabilities for two different groups, as compared with Laplace´s method for one homogeneous group of the same number of observations.

The findings, displayed in Tab. 1, contribute to the question, why variety of evidence permits more reliable conclusions than the same number of Es of a homogeneous semantic content (For literature see e.g. J. Earman, 1992, Ch 3 and passim). It may be assumed that the difference between the results of the Laplacean method and of the cMPE-Method will be even more marked, if the observations, pertaining to the latter group will be partitioned into more than two semantically diverse groups.

IV: FUNDAMENTAL DIFFERENCES BETWEEN THE BAYESIAN METHOD AND THE cMPE-METHOD.

These differences may be also noted in the Figs. 2-4.

*A: Conditional and non-conditional probabilities.*

At least one premise of BayT [the likelihood, P(E|H)] must be a conditional probability and conditional probabilities are found in the numerator, as well as in the denominator of (2a and b),. The cMPE-Method may work without conditional (i.e. with a-priori) probabilities.

B: *Divisions and Additions*
In BayT an increase in probability is achieved by division of data by a sum of products, <1.
By the cMPE-Method two non-disjunctive probabilities are added, however, taking into account their non-linearity.

C: *Ascending and descending operations.*
As shown in Fig. 2 **d**, by BayT descending operations (the result smaller than the prior probabilities) are possible. This will be frequently the case if P(B|¬A) is large and P(A) x P(B|A) is small. By the cMPE-Method only ascending operations are possible. Both of the operations may be applied if A supports B or B supports A actively, however, they must not be applied for coincidence by sheer chance.

D: *Dependences.*
In BayT at least one conditional probability enters the operation and the result again is a conditional probability. Thus, no stochastic independence is required and possible. Also for the cMPE-Method no stochastic independence is required, however, semantic independence, as described above, is prerequisite for the multiplicative procedure.

E: *Probabilities 1 or nearly 1.*
Most of the philosophers and natural scientists adhere to the Humean tradition by denying our capacity to arrive inductively at universal propositions or laws of nature. By BayT, however, you easily arrive at probabilities of 1.0, and the same is true for the principle of relative frequency. There just may be no counterexample (no P(E|¬H) and we arrive after a few observations by (2a or b) at 1.0. Similarly, in relative frequency, a few observations of just favourable cases may result in P = 1.0 . Of course, if high probabilities result from small numbers of observation, we always will be extremely suspicious.

This is different for Laplace´s method [P = (n + 1) / (n + 2)] and for the cMPE-Method. In both of the methods, the probabilities achieved approach probability 1 asymptotically, whereby the Humean scruples are complied with. No probability 1.0 may be arrived at by the cMPE-Method, unless you add 1.0 .

V THE ROLE OF WIGHT BEARING AND OF EXTENSIONAL EVIDENCE
A *General Considerations.*
Some closer look at the different roles, some evidence can play, may improve our understanding of probabilistic processes. Here, first a comparison of some probability formulas.

| $P = \dfrac{N^F}{N^P}$ | $P = \dfrac{m+1}{n+2}$ | $P(H\|E) = \dfrac{P(H) \times P(E\|H)}{P(H) \times P(E\|H) + P(\neg H) \times P(E\|\neg H)}$ | $P(H\|E) = \dfrac{P(H)}{P(E)}$ |
|---|---|---|---|
| **a** | **b** | **c** | **d** |

Fig. 5: a: : The method of relative frequency. $N^F$ = number of favourable cases, $N^P$: Number of possible cases. b: Laplace`s Formula. m: The number of white balls drawn from an urn; n: The total

number of drawn balls. c: The BayT formula. d: The BayT formula, however for conditions of H implying E.

In Fig. 5 we see from a - d always the same principle: In the numerator of a fraction the weight-bearing E, that is favourable to our interests, and in the denominator the extension, the number of possible cases, among which the favourable cases are dispersed. Also in c, the BayT-formula, this principle is preserved, although the numerator and the denominator are expanded. According to (1b) we even may write this formula: $P(H|E) = P(H \& E) / P(H \& E) + P(\neg H \& E)$, whereby this formula resembles the relative frequency formula even more. Typically, the numerator of the BayT formula contains the weight bearing numbers in form of a multiplication of $P(H)$ and $P(E|H)$ as factors, whereby the product indicates the overlapping area in Fig. 2 c and d. The denominator contains the extension in form of an addition. This principle is simplified in d, because due to the implication, $P(E|H)$ became 1. However, there is a fundamental difference between formulas a or b on the one hand, and c and d on the other: a and b have numbers >1 in the denominator and therefore, by the division, the whole fraction is diminished, whereas in c and d, because of numbers <1 in the denominator, the value of the whole fraction increases. Some more pitfalls are connected with the implication (d) and these will be studied in the next subsection.

B: *Implication of E by H.*

Two different constellations may be encountered:

(i): *Hypothesis H implies one single evidence.* Cases with $P(E|H) = 1$ are no rarities. They are generally calculated according to (8):

$$P(H|E) = P(H) / P(E) \qquad (8)$$

As in Fig. 5, we have on the right-hand side of (8) the weight-bearing number in the numerator and the extension in the denominator, and this is correct, as long as $P(E)$ really indicates an extension. However, quite frequently in scientific research, $P(E)$ may indicate the weight-bearing evidence, as shown by the following

*Example:* It is rumoured that René Descartes 1650 died at the court of Kristina of Sweden due to being poisoned by arsenic, the standard poison in renaissance age. Arsenic was also assumed because his death was preceded by severe enteritis2 and as a rule, after arsenic poisoning severe enteritis is observed. Here, we are interested only in the probabilistic handling. Two different possibilities are indicated in Fig. 6

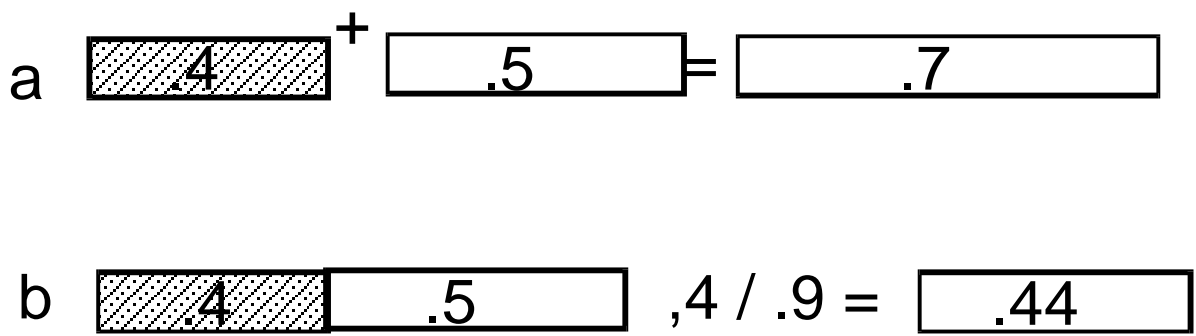

Fig. 6: Implication of E by H. a: Weight-bearing approach. P(E), the observation of severe enteritis supports P(H), arsenic poisoning. b: Extensional approach: $P(H|E) = P(H) / P(E_{total})$ according to (8).

In Fig. 6a we assume that the occurrence of enteritis supports the probability of arsenic poisoning. We first allot .4 to the hypothesis P(H) and ,5 to the weight-bearing data P($E_W$) and arrive by the cMPE-Method at .7. However, in Fig. 6b, we see the BayT-approach. We assume that besides arsenic poisoning, there might be other causes of enteritis (.5) [2] and arrive at P(H) / [P(H) + P(E|¬H)] = .4 / .9 = .44 according to (8). This is a comparatively small support that does not allow for the active support, as may be given by the weight-bearing E. On the other hand, high probabilities may be arrived at according to (9), if P(E|¬H) were small or zero. Thus we may say that the active support according to the cMPE-Method (Fig. 6a) neglects the dilution of the probabilities by P($E_{EXT}$,), especially by P(E|¬H) and the BayT according to (8) (Fig, 6b) neglects partially the active support, taking into account only P(H) and P($E_{EXT}$), but neglecting the weight-bearing E.

How can we avoid these shortcomings? What we need, is a procedure that allows for the weight-bearing evidence P($E_W$) as well as for the extensional evidence P($E_{EXT}$). Such a method will be described in a forthcoming paper.

*Weight-bearing probabilities* are all probabilities of evidence that positively contribute to the value of a fraction. To them belong causally determined probabilities, such as the incidence of thunderstorms in dependence of air humidity or the incidence of enteritis after arsenic poisoning. Also probabilities arrived at by the method of relative frequency or by Laplace´s method are weight-bearing. The latter operations are straightforward converters of extensional probabilities (found in their denominator) into weight- bearing probabilities. On the other hand, probabilities gained by conventional addition of a multitude of probabilities may become extensional as in the denominator of BayT (2a and b).

Only weight-bearing probabilities may enter into cMPE-calculations. Here, P(E) doesn't show up in the denominator of an equation but enters as complement into a product of complements. Therefore, a higher weight-bearing P(E) increases p(H) $^+$ P(E) more than a small one, however, these data cannot be applied for (9), but must be calculated by the cMPE-Method. Originally, a Bayesian solution of this question was tried, however, in the course of these investigations, quite general doubts concerning the Bayesian algorithm arose. Therefore, this part of the paper was cancelled. The present amended paper may serve as an introduction to a forthcoming paper in which I shall derive the rMPE formula, a new procedure for increasing the probabilities of hypotheses by evidence.

(ii) *Hypothesis H implies a multitude of evidences.* C.D. Broad already in 1918 commented on this question. He wrote:

*"We see that the probability of a hypothesis is increased as we verify its consequences because the initial probability is the numerator of a fraction whose denominator is a product which contains more factors (and since they are proper fractions, grows  s m a l le r) the more consequences we deduce and verify".* (Original emphasis).

In what follows, I copy two formulas of C.D. Broad´s that presuppose implication (the second is an expansion of the denominator of the first):

$h|c_1 c_2 \ldots c_n|f = h|f / c_1 c_2 \ldots c_n|f$

$c_1 c_2 \ldots c_n|f = c_1|f \times c_2|c_1 f \times c_3|c_2 c_1 f \times \ldots \times c_n|c_{n-1} \ldots c_1 f$

---

[2] Jesuits trying to convert the queen to Catholicism. were blamed for this crime. According to other reports, Descartes died because of pneumonia, caused by his obligation to lecture every day at 5 o clock a.m.

(h: hypothesis; cn: probability of the verified prediction n; f: Background knowledge) Here three questions must be answered: 1. Must we tolerate that if we verify numerous consequences of the hypothesis, the denominator eventually will become smaller than the numerator and the whole fraction will become >1?[3] 2. What rule permits us to multiply (instead of adding like in 2a or in 2b) the data in the denominator, whereby the denominator decreases? And 3. How can we explain that data that must be regarded as weight-bearing probabilities appear in the denominator?-Our qualms will be illustrated by

*anotherExample:*(Gottlob, 2000) There, I considered the discovery of certain bacteria as causes of tuberculosis by Robert Koch. In the analysis I assumed a prior probability of the hypothesis of .5. Koch investigated on four different predictions and was able to verify them, mostly experimentally. I allotted these predictions a comparatively high probability of .8 each. In spite of this high rating we arrive by the above Broad´s formula at .5 / .8 x.8 x.8 x.8 = 1.2

As the main cause for this inconsistency I see the accumulation of (weight-bearing) probabilities in the denominator without increasing adequately the weight-bearing probabilities in the numerator. Also for these problems the forthcoming paper will offer a solution.

VI: CONCLUSIONS

A: There is no competition between the Bayesian Method and the cMPE-Method. Whereas for the BayT conditional probabilities are prerequisite and are also arrived at, the cMPEMethod may do with non-conditional (a-priori) probabilities.

B Some shortcomings of the BayT:

(i) Parts of a hypothesis that cannot enter into the BayT-operation (see `R´ in Fig. 2c) may be added by the cMPE-Method to the result.

(ii) If in BayT operations hypothesis H implies the evidence E, and evidence E is weight bearing, the conventional operation according to (8) will yield wrong results.

(iii) If Hypothesis (H) implies a multitude of predictions ($E_1$, $E_2$, $E_n$) that can be verified, the conventional procedure of multiplying the evidences in the denominator of a fraction may yield probabilities >1 and is therefore flawed.

C: The Old Evidence Puzzle arises due to the non-linearity of non-disjunctive probabilities. High probabilities, added by high probabilities, underlie but small changes if judged according to a linear scale. However, the real increase in probability may be evaluated by the non-linear DPE-Method.

D: There are many possibilities of application of the cMPE- or the DPE-Method. Just briefly mentioned were forensic applications or the evaluation of the `therapeutic window´ in medicine. The MPE-Method is likely to underlie many, mostly unconscious cognitive processes inclusive Gestalt perception. I close with a hint at a possible philosophical significance: Every human being may notice in any second of its waking life that all material objects underlie the gravitational attraction. Exceptions, objects, flying in the air, are easily explained by the laws of buoyancy or by the aerodynamic lift. Six billions of living human beings of an average age of say, 25 years have each, at least observed a million of objects that were lying on their support. Thus we arrive at least at $6 \times 10^{15}$ observations without any unexplained for counter instance and at a probability of error less than $10^{-15}$. A similar number of observations may be obtained for the laws of conservation and for the uniformity of nature. Again for both of these principles, we may assume probabilities of error of at most $10^{-15}$. The three laws mentioned (and several more) support each other. Thus we may multiply the probabilities of error and arrive at a common probability of error of at most $10^{-45}$. I venture that here may be hidden an answer to Hume´s problem. Hume is right in doubting our ability to arrive at certain universal propositions and laws of nature. However, in his time, he was unable to realize the mighty power behind probabilistic reasoning that is based on overwhelming numbers of observation and that permits reducing probabilities of error to an extreme minimum.

---

[3] To a similar denominator, different only by a different notation, J. Earman (1992,p.107) writes:
If $Pr(H|K) > 0$, the denominator on the right-hand side will eventually become smaller than the numerator, which contradicts an axiom of probability, unless $Pr(E_{n+1}|\&_{i \leq n} E_i \& K) \to 1$ as $n \to \infty$.

*Acknowledgement:* The author feels greatly indebted to Prof. Viktor Scheiber, former head of the Institute of Medical Statistics, Vienna University.

Author´s E-Mail: rainer.gottlob@univie.ac.at

REFERENCES
Boole, G. (1857) On the Application of the Theoriy of Probabilities to the Question of the Combination of Testimonies or Judgements. *Transactions of the Royal Society of Edinburgh, XXI,* quoted by L.J. Cohen (1976)
Boole, G. (1953) *Studies in Logic and Probability*, London, Watts & Co (1976)
Broad, C.D. (1918) On the Relation between Induction and Probability. *Mind (NS)* 108, pp. 388 - 404
Cohen, L. J. (1977) *The Probable and the Provable*. Clarendonn Prss, Oxford
Cohen, L. J. (1976) How can one Testimony Corroborate Another? In *Essays in Memory of Imre Lakatos* (ed.: Cohen, R.S., Feyedrabend, P.-K., Wartofsky, M.W.), pp. 65 78 . Dordrecht, Reidel,
Earman, J. (1992) *Bayes or Bust?* Bradford Book, MIT-Press, Cambridge, Mass., London.
Ekelöf, P. O. (1964) Free Evaluation of Evidence. *Scandinavian Studies in Law*, 8, 47 65
Ekelöf, P.O. (1981) Beweiswert in *Festschrift für Franz Baur* (ed. Grunsky, W., Stürner, R., Walter, G. and Wolf, M) pp. 343 363, Tübingen, Mohr
Glymour, C. (1980) *Theory and Evidence*. (Princeton, Princeton Univ. Press
Gottlob, R. (2000) New Aspects of the Probabilistic Evaluation of Hypotheses and Experience. *International Studies in the Philosophy of Science* 14., 147 - 163
Gottlob, R., (2003) Algorithm for Addition and Subtraction of Two or More Non conditional Probabilities by the MPE or the DPE-Method and for Calculating Probabilities by Percentages. *www. mathpreprints.com/math/Preprint/Gottlob/* 2003.2/3
Gottlob, R. (2004) Die Rolle einer multiplikativen Wahrscheinlichkeitsmethode b ei der Gestaltwahrnehmung. In *Experience and Analysis, Papers of the 27$^{th}$ International Wittgenstein Symposium* , E. Leinfellner, R. Haller, W. Leinfellner, K. Puhl and P. Weingartner (eds.) Vol. XII, pp. 129-121, Kirchberg/Wechsel
Keynes, J. M. (1957) *A Treatise on Probability*. Macmillan & Co, London
Mises, von, R. (1990) *Kleines Lehrbuch des Positivismus* (R. Stadler, (ed.) Suhrkamp, Frankfurt/Main
Russell, B (1948) *Human Knowledge, its Scope and Limits*. Allen & Unwin, London.